\documentclass{amsart}
\usepackage{latexsym}
\usepackage{amsmath,amssymb,stmaryrd}
\usepackage{amscd,amsxtra}
\usepackage{amsthm}

\begin{document}
\title{On the highest Lyubeznik number of a local ring}
\author{Wenliang Zhang}
\address{Department of Mathematics, University of Minnesota, Minneapolis, MN 55455}
\email{wlzhang@math.umn.edu}
\maketitle

\begin{abstract}
Let $A$ be a $d$-dimensional local ring containing a field. We will prove that the highest Lyubeznik number $\lambda_{d,d}(A)$ (defined in \cite{l1}) is equal to the number of connected components of the Hochster-Huneke graph (defined in \cite{hh}) associated to $B$, where $B=\widehat{\hat{A}^{sh}}$ is the completion of the strict Henselization of the completion of $A$. This was proven by Lyubeznik in characteristic $p>0$. Our statement and proof are characteristic-free.
\end{abstract}

\section{Introduction}
Throughout this paper, all rings are Noetherian and commutative. Let $A$ be a local ring that admits a surjection from an $n$-dimensional regular local ring $(R,\mathfrak{m})$ containing a field. Let $I\subset R$ be the kernel of the surjection, and let $k=R/{\mathfrak{m}}$ be the residue field of $R$. Then the Lyubeznik numbers $\lambda_{i,j}(A)$ (Definition 4.1 in \cite{l1}) are defined to be $dim_k(Ext_R^i(k,H^{n-j}_I(R)))$. And it was proven in \cite{l1} that they are all finite and depend only on $A,i$ and $j$, but neither on $R$, nor on the surjection $R\rightarrow A$. \par
The Lyubeznik numbers have been studied by a number of authors, including \cite{k1}, \cite{k2}, \cite{l3}, \cite{ls}, \cite{w}. In this paper, we will give an interpretation of $\lambda_{d,d}(A)$ in terms of the topology of $SpecA$.\par
Firstly, we reproduce the definition of the Hochster-Huneke graph associated to a local ring which was originally given in \cite{hh}:
\newtheorem*{def1}{{\bf Definition 1.1} (Definition 3.4 in \cite{hh})}
\begin{def1}
Let $B$ be a local ring. The graph $\Gamma_B$ associated to $B$ is defined as follows. Its vertices are the top-dimensional  minimal prime ideals of $B$, and two distinct  vertices $P$ and $Q$ are joined by an edge if and only if $ht_B(P+Q)=1$. 
\end{def1}

In \cite{l2} and \cite{l3}, the following question was posed
\newtheorem*{q2}{{\bf Question 1.2} (Question 1.1 in \cite{l3})}
\begin{q2}
Is $\lambda_{d,d}(A)$ equal to the number of the connected components of the Hochster-Huneke graph $\Gamma_B$ associated to 
$B=\widehat{\hat{A}^{sh}}$, the completion of the strict Henselization of the completion of $A$?
\end{q2} 

As is pointed out in \cite{l3}, the graph $\Gamma_B$ can be realized by a much smaller ring than $B=\widehat{\hat{A}^{sh}}$. Namely, if $\hat{A}$ is the completion of $A$ with respect to the maximal ideal and $k\subset \hat{A}$ is a coefficient field, then there exists a finite separable extension field $K$ of $k$ such that $\Gamma_B=\Gamma_{\hat{A}\otimes_kK}$. In particular, if the residue field of $A$ is separably closed, then $\Gamma_B=\Gamma_{\hat{A}}$.\par

It is shown in \cite{l3} that the answer to the above question is positive in characteristic $p>0$.
Our main result in this paper is that the answer to the above question is positive in general, without any restriction on the characteristic, i.e,

\newtheorem*{thm}{\bf Main Theorem}
\begin{thm}
Let $A=R/I$ be a local ring, where $R$ is a regular local ring containing a field (of any characteristic), and $dim(A)=d$. Then $\lambda_{d,d}(A)$ is equal to the number of connected components of the Hochster-Huneke graph $\Gamma_B$ associated to $B=\widehat{\hat{A}^{sh}}$.
\end{thm}

Our proof of the Main Theorem is completely characteristic-free. We use the following result from \cite{l3} whose proof in \cite{l3} is completely characteristic-free.
\newtheorem*{lem1.3}{{\bf Lemma 1.3} (Corollary 2.4 in \cite{l3})}
\begin{lem1.3}
Let $A$ be a local ring of dimension $d$ containing a field and let $B=\widehat{\hat{A}^{sh}}$ be the completion of the strict Henselization of the completion of $A$. Let $\Gamma_1,\cdots, \Gamma_r$ be the connected components of $\Gamma_B$. Let $I_j$ be the intersection of the minimal primes of $B$ that are the vertices of $\Gamma_i$. Let $B_j=B/I_j$. Then $\lambda_{i,d}(A)=\sum^{r}_{j=1}\lambda_{i,d}(B_j)$ for every $i$.
\end{lem1.3}

Clearly, to prove our Main Theorem it is enough to show that $\lambda_{d,d}(B_j)=1$ for every $j$. Since every $B_j$ is complete, local, $d$-dimensional, reduced, equidimensional, contains a field, has a separably closed residue field, and the Hochster-Huneke graph associated to $B_j$ is connected, this is proven in the following Theorem 1.4.

\newtheorem*{thm2}{\bf Theorem 1.4}
\begin{thm2}
If $A$ is $d$-dimensional, complete, reduced, equidimensional, local, contains a field, has a separably closed residue field and the graph associated to $A$ is connected, then $\lambda_{d,d}(A)=1$.
\end{thm2}

Thus our Main Theorem follows from Lemma 1.3 and Theorem 1.4. To complete the proof of the Main Theorem it remains to prove Theorem 1.4. This is accomplished in the following section.\par
\section*{Acknowledgements}
The result of this paper is from my thesis. I would like to thank my advisor Prof. Gennady Lyubeznik for suggesting this problem to me.
\section{Proof of Theorem 1.4}
Throughout this section $A$ is as in Theorem 1.4, i.e. $d$-dimensional, complete, local, reduced, equidimensional, contains a field, has a separably closed residue field and the Hochster-Huneke graph associated to $A$ is connected. The case that $dim(A)\leq 2$ has been completely settled by Kawasaki \cite{k1} and Walther \cite{w}, independently. Thus it remains to settle the case when $dim(A)\geq 3$. We will do this by induction on $dim(A)$, the case that $dim(A)\leq 2$ being known.\par
Accordingly, thoughout this section we assume that $d\geq 3$ and Theorem 1.4 proven for $d-1$. By Cohen's Structure Theorem, $A$ is a homomorphic image of a complete regular local ring ($R,\mathfrak{m}$) containing a field.\par

Let $dim(R)=n$. By \cite[3.6]{l1} the set of the minimal primes of the support of $H^{n-d+1}_I(R)$ is finite. Hence standard prime avoidance implies that there is an element $r\in\mathfrak{m}$ that does not belong to any minimal prime of $I$ nor to any minimal prime of the support of $H^{n-d+1}_I(R)$ different from $\{\mathfrak{m}\}$ (if $\mathfrak{m}$ is a minimal (hence the only associated) prime of $H^{n-d+1}_I(R)$, then the only condition on $r$ is that $r\in\mathfrak{m}$ and $r$ does not belong to any minimal prime of $I$). We fix one such element $r\in\mathfrak{m}$ throughout the rest of this section.\par

Let $\bar{r}$ be the image of $r$ in $A=R/I$. Then $\bar{r}$ is not contained in any minimal prime ideal of $A$, since $r$ is not contained in any minimal prime of $I$. Hence $A/\bar{r}A$ is equidimensional and $dim(A/\bar{r}A)=d-1$. We are going to prove the following two propositions.

\newtheorem*{p2.1}{\bf Proposition 2.1}
\begin{p2.1}
$\lambda_{d,d}(A)=\lambda_{d-1,d-1}(A/\sqrt{\bar{r}A})$.
\end{p2.1}
\newtheorem*{p2.2}{\bf Proposition 2.2}
\begin{p2.2}
The Hochster-Huneke graph associated to $A/\sqrt{\bar{r}A}$ is connected.
\end{p2.2}
\begin{proof}[Proof of Theorem 1.4]
Proposition 2.2 shows that the ring $A/\sqrt{\bar{r}A}$ satisfies the hypotheses of Theorem 1.4 for dimension $d-1$. The inductive hypotheses (i.e. Theorem 1.4 in dimension $d-1$) implies that $\lambda_{d-1,d-1}(A/\sqrt{\bar{r}A})=1$. Proposition 2.1 now shows that $\lambda_{d,d}(A)=1$. This completes the proof of Theorem 1.4 modulo Propositions 2.1 and 2.2.
\end{proof}
It remains to prove Propositions 2.1 and 2.2. We begin with a proof of Proposition 2.1 which requires Lemmas 2.3 and 2.4 below.

\newtheorem*{l2.3}{\bf Lemma 2.3}
\begin{l2.3}
$dim(Supp_R(H^{n-d+1}_I(R)))\leq d-2.$
\end{l2.3}
\begin{proof}[Proof]
Let $P$ be an arbitrary element of $Supp_R(H^{n-d+1}_I(R))$. Assume $ht(P)\leq n-d+1$. If $ht(P)\leq n-d$, then $(H^{n-d+1}_I(R))_P\cong H^{n-d+1}_{IR_P}(R_P)=0$ since the dimension of $R_P$ is less than $n-d+1$. If $ht(P)=n-d+1$, then the dimension of $R_P$ is $n-d+1$ and by the Hartshorne-Lichtenbaum Vanishing Theorem (\cite[8.2.1]{bs}) $$(H^{n-d+1}_I(R))_P=H^{n-d+1}_{IR_P}(R_P)=0$$
(since $I\widehat{R_P}$ is not $P\widehat{R_P}$-primary as $R$ is regular and every minimal prime of $IR_P$ has height $n-d$). So, $ht(P)\geq n-d+2$ for every $P\in Supp_R(H^{n-d+1}_I(R))$, hence $dim(Supp_R(H^{n-d+1}_I(R)))\leq d-2$.
\end{proof}

\newtheorem*{l2.4}{\bf Lemma 2.4}
\begin{l2.4}
$dim(Supp_R(H^0_{(r)}(H^{n-d+1}_I(R))))\leq d-3.$
\end{l2.4}
\begin{proof}[Proof]
If $Supp_R(H^{n-d+1}_I(R))=\{\mathfrak{m}\}$, then $H^{n-d+1}_I(R)$ is injective by \cite[3.6]{l1}. Then $dim(Supp_R(H^0_{(r)}(H^{n-d+1}_I(R))))=0\leq d-3$, since $d\geq 3$ by our assumption. If $Supp_R(H^{n-d+1}_I(R))\neq \{\mathfrak{m}\}$, let $P$ be an arbitrary element of $Supp_R(H^0_{(r)}(H^{n-d+1}_I(R)))$. Then $P$ has to contain $r$ and a minimal element $P'$ in $Supp_R(H^{n-d+1}_I(R))$. Since $r\notin P'$ by Lemma 2.1, it follows from Krull's Principal Ideal Theorem that $ht(P'+(r))=ht(P')+1$. Hence, $ht(P)\geq n-d+2+1=n-d+3$. Therefore, $dim(Supp_R(H^0_{(r)}(H^{n-d+1}_I(R))))\leq d-3$.
\end{proof}

\begin{proof}[Proof of Proposition 2.1]
We have the Grothendieck spectral sequence for the composition of functors 
$$E^{p,q}_2 =H^p_{(r)}(H^q_I(R))\Rightarrow H^{p+q}_{I+(r)}(R).$$
In this spectral sequence, all differentials $d^{p,q}_s:E^{p,q}_s\rightarrow E^{p+s,q-(s-1)}_s$ are zero since $H^{p+s}_{(r)}(H^{q-(s-1)}_I(R))$ is zero (indeed $s\geq 2$ implies $p+s\geq 2$ and the ideal ($r$) is 1-generated). Therefore
$$E^{p,q}_{\infty}=E^{p,q}_2=H^p_{(r)}(H^q_I(R)).$$
The spectral sequence is convergent, hence we have a finite filtration
$$0=F^tH^m_{I+(r)}(R)\subseteq \cdots F^{p+1}H^m_{I+(r)}(R)\subseteq F^pH^m_{I+(r)}(R)\cdots \subseteq F^sH^m_{I+(r)}(R)=H^m_{I+(r)}(R)$$
so that 
$$E^{p,q}_{\infty}\cong F^pH^{p+q}_{I+(r)}(R)/F^{p+1}H^{p+q}_{I+(r)}(R).$$
Since $E^{p,q}_2=0$ for any $p>1$, we have 
$$F^2H^m_{I+(r)}(R)=\cdots =F^tH^m_{I+(r)}(R)=0.$$
Therefore the filtration is in fact
$$0=F^2H^m_{I+(r)}(R)\subseteq F^1H^m_{I+(r)}(R)\subseteq F^0H^m_{I+(r)}(R)=H^m_{I+(r)}(R),$$
and 
$$E^{0,m}_2\cong F^0H^m_{I+(r)}(R)/F^1H^m_{I+(r)}(R),$$
$$E^{1,m-1}_2\cong F^1H^m_{I+(r)}(R).$$
In particular, we have an exact sequence
$$\ \ 0\rightarrow E^{1,m-1}_2\rightarrow H^m_{I+(r)}(R)\rightarrow E^{0,m}_2\rightarrow 0.$$
i.e. we have the following exact sequence
$$(*)\qquad 0\rightarrow H^1_{(r)}(H^{m-1}_I(R))\rightarrow H^m_{I+(r)}(R)\rightarrow H^0_{(r)}(H^m_I(R))\rightarrow 0.$$
Applying $\Gamma_{\mathfrak{m}}$ to ($*$) and setting $m=n-d+1$, we have a long exact sequence 
$$\cdots\rightarrow H^{d-2}_{\mathfrak{m}}(H^0_{(r)}(H^{n-d+1}_I(R)))
\rightarrow H^{d-1}_{\mathfrak{m}}(H^1_{(r)}(H^{n-d}_I(R)))
\rightarrow$$ 
$$H^{d-1}_{\mathfrak{m}}(H^{n-d+1}_{I+(r)}(R))\rightarrow H^{d-1}_{\mathfrak{m}}(H^0_{(r)}(H^{n-d+1}_I(R)))\rightarrow\cdots.$$
Since $dim(Supp_R(H^0_{(r)}(H^{n-d+1}_I(R))))\leq d-3$,
$$ H^{d-2}_{\mathfrak{m}}(H^0_{(r)}(H^{n-d+1}_I(R)))=H^{d-1}_{\mathfrak{m}}(H^0_{(r)}(H^{n-d+1}_I(R)))=0,$$
by Grothendieck's Vanishing Theorem \cite[6.1.2]{bs}.  Hence, the above long exact sequence implies
$$(**)\qquad H^{d-1}_{\mathfrak{m}}(H^1_{(r)}(H^{n-d}_I(R)))\cong H^{d-1}_{\mathfrak{m}}(H^{n-d+1}_{I+(r)}(R)).$$
The height of $I+(r)$ is $n-d+1$ and $R$ is regular, thus 
$$H^{n-d}_{I+(r)}(R)=0.$$
So, if we set $m=n-d$ in ($*$), we will have 
$$H^0_{(r)}(H^{n-d}_I(R))=0.$$
For any $R$-module, we always have the following exact sequence
$$(***)\qquad 0\rightarrow H^0_{(r)}(M)\rightarrow M\rightarrow M_r\rightarrow H^1_{(r)}(M)\rightarrow 0.$$
Let $M=H^{n-d}_I(R)$ in ($***$), then we will have an exact sequence
$$(****)\qquad 0\rightarrow H^{n-d}_I(R)\rightarrow (H^{n-d}_I(R))_r\rightarrow H^1_{(r)}(H^{n-d}_I(R))\rightarrow 0.$$
Applying $\Gamma_{\mathfrak{m}}$ to ($****$) and keeping in mind that $H^j_{\mathfrak{m}}((H^{n-d}_I(R))_r)=0$ for $j>0$ (since mutliplication by $r\in \mathfrak{m}$ is bijective on $H^{n-d}_I(R)_r$), we have
$$H^{d-1}_{\mathfrak{m}}(H^1_{(r)}(H^{n-d}_I(R)))\cong H^d_{\mathfrak{m}}(H^{n-d}_I(R)),$$ 
as $d\geq 3$. Considering ($**$), we have
$$H^d_{\mathfrak{m}}(H^{n-d}_I(R))\cong H^{d-1}_{\mathfrak{m}}(H^{n-d+1}_{I+(r)}(R)).$$
Since $H^{n-d+1}_{I+(r)}(R)=H^{n-d+1}_{\sqrt{I+(r)}}(R)$, we have 
$$dim_{R/\mathfrak{m}}Hom_R(R/\mathfrak{m},H^d_{\mathfrak{m}}(H^{n-d}_I(R)))=dim_{R/\mathfrak{m}}Hom_R(R/\mathfrak{m},H^{d-1}_{\mathfrak{m}}(H^{n-(d-1)}_{\sqrt{I+(r)}}(R))),$$
which implies that $$\lambda_{d,d}(A)=\lambda_{d,d}(R/I)=\lambda_{d-1,d-1}(R/(\sqrt{I+(r)})=\lambda_{d-1,d-1}(A/\sqrt{\bar{r}A}).$$ This completes the proof of Proposition 2.1.
\end{proof}

To complete the proof of Theorem 1.4, it remains to prove Proposition 2.2. Let $\Theta=\{P_1,\cdots,P_s\}$ be the set of the minimal prime ideals of $A$. $\text{Let}\ \Sigma_i =\{Q\in Spec(A)|Q\ is\ minimal\ over\ P_i+\sqrt{\bar{r}A}\}$, and let $\Sigma=\cup_i\Sigma_i$. There is 1-1 correspondence between $\Sigma$ and the set of the minimal prime ideals in $A/\sqrt{\bar{r}A}$.\par
We recall the the height of an ideal is the minimum of the heights of the minimal primes over that ideal.
\newtheorem*{l2.5}{\bf Lemma 2.5}
\begin{l2.5}
Let $P_1$ and $P_2$ be two arbitrary elements in $\Theta$. If  for any $Q_{\alpha}\in \Sigma_1$ and any $Q_{\beta}\in \Sigma_2$, $ht_{A/\sqrt{\bar{r}A}}((Q_{\alpha}+Q_{\beta})/\sqrt{\bar{r}A})\geq 2$, then $ht_{A}(P_1+P_2)\geq 2.$
\end{l2.5}
\begin{proof}[Proof]
Otherwise, $ht_{A}(P_1+P_2)=1$ (obviously, $ht_{A}(P_1+P_2)\geq 1$). By the Principal Ideal Theorem and considering that $A$ is catenary because it is complete, we have 
$$(\bullet)\qquad ht_{A}(P_1+P_2+\sqrt{\bar{r}A})\leq 2.$$
Let $\mathcal{Q}$ be an arbitrary prime ideal of $A$ minimal over $P_1+P_2+\sqrt{\bar{r}A}$. Then $\mathcal{Q}$ must contain some $Q_1\in \Sigma_1$ and some $Q_2\in \Sigma_2$. Therefore,
$$ht_{A/\sqrt{\bar{r}A}}(\mathcal{Q}/\sqrt{\bar{r}A})\geq ht_{A/\sqrt{\bar{r}A}}((Q_1+Q_2)/\sqrt{\bar{r}A})\geq 2.$$
Thus,
$$ht_{A/\sqrt{\bar{r}A}}((P_1+P_2+\sqrt{\bar{r}A})/\sqrt{\bar{r}A})\geq 2.$$
Hence, for any prime ideal $\tilde{Q}$ minimal over $P_1+P_2+\sqrt{\bar{r}A}$, there exist prime ideals $Q$ and $\bar{Q}$ so that we have a chain of ideals 
$$\sqrt{\bar{r}A}\subset Q\subsetneq \bar{Q}\subsetneq \tilde{Q}.$$
$Q$ contains $\sqrt{\bar{r}A}$, thus $Q$ properly contains some element in $\Theta$, say, $P_3$. Then we have a chain of ideals 
$$P_3\subsetneq Q\subsetneq \bar{Q}\subsetneq \tilde{Q}.$$
Therefore, for every prime ideal $\tilde{Q}$ minimal over $P_1+P_2+\sqrt{\bar{r}A}$,
$$ht_{A}(\tilde{Q}) \geq 3,$$
hence,
$$ht_{A}(P_1+P_2+\sqrt{\bar{r}A})\geq 3,$$
contrary to ($\bullet$).
\end{proof}

\newtheorem*{l2.6}{\bf Lemma 2.6}
\begin{l2.6}
Assume that the Hochster-Huneke graph associated to $A/\sqrt{P_i+\bar{r}A}$ is connected for all $P_i\in\Theta$. Then so is the Hochster-Huneke graph associated to $A/\sqrt{\bar{r}A}$.
\end{l2.6}
\begin{proof}[Proof]
Indeed, assume the graph associated to $A/\sqrt{\bar{r}A}$ is not connected. Then $\Sigma$ can be divided into 2 non-empty disjoint subsets: $\Sigma^*$ and $\Sigma^{**}$, such that if $Q_1\in \Sigma^*$ and $Q_2\in\Sigma^{**}$ then $ht_{A/\sqrt{\bar{r}A}}((Q_1+Q_2)/\sqrt{\bar{r}A})\geq 2$. Since the graph associated to $A/\sqrt{P_i+\bar{r}A}$ is connected, $\Sigma_i$ has to be completely contained in $\Sigma^*$ or completely contained in $\Sigma^{**}$ for every $i$. Therefore, we can divide $\Theta$ into 2 non-empty disjoint subtsets:
$$\Theta^*=\{P\in\Theta|the\ prime\ ideals\ minimal\ over\ P+\sqrt{\bar{r}A}\ are\ contained\ in\ \Sigma^* \}$$ 
and
$$\Theta^{**}=\{P\in\Theta|the\ prime\ ideals\ minimal\ over\ P+\sqrt{\bar{r}A}\ are\ contained\ in\ \Sigma^{**} \}.$$ 
For an arbitrary element $P_1\in \Theta^*$ and an arbitrary element $P_2\in \Theta^{**}$, $$ht_{A/\sqrt{\bar{r}A}}((Q_1+Q_2)/\sqrt{\bar{r}A})\geq 2,\ \forall Q_1\in \Sigma_1,\forall Q_2\in \Sigma_2.$$
Lemma 2.5 implies $ht_A(P_1+P_2)\geq 2$. Hence, the graph associated to $A$ is not connected either, which is a contradiction.
\end{proof}

\begin{proof}[Proof of Proposition 2.2]
According to Lemma 2.6, it is enough to prove that the graph associated to $A/\sqrt{P_i+\bar{r}A}$ is connected for all $P_i\in \Theta$. Denoting $A/P_i$ by $A$ and the image of $\bar{r}$ in $A/P_i$ by $\bar{r}$ again, we are reduced to proving that if $A$ is a domain and $\bar{r}\in \mathfrak{m}$ is nonzero, then the graph associated to $A/\sqrt{\bar{r}A}$ is connected.\par
The following result is not explicitly stated in \cite{hh}, but is a straightforward consequence of \cite[3.6c,e]{hh} and \cite[3.9b,c]{hh}: \emph{Let $S$ be a complete local equidimensional ring. If $S$ satisfies Serre's condition $S_2$ and $x_1,\cdots,x_k$ is a part of a system of parameters of $S$, then the graph associated to $S/\sqrt{(x_1,\cdots,x_k)}$ is connected.}\par
Let $S$ be the normalization of $A$. Since $A$ is a complete local domain, so is $S$. Serre's criterion of normality shows that $S$ is $S_2$. Since $\bar{r}\in S$ is $S$-regular, the graph associated to $S/\sqrt{\bar{r}S}$ is connected by the above-quoted result. As $S$ is module-finite over $A$, the going-up theorem implies that $\sqrt{\bar{r}S}\cap A=\sqrt{\bar{r}A}$, hence the natural map $\phi:A/\sqrt{\bar{r}A}\rightarrow S/\sqrt{\bar{r}S}$ is injective and $S/\sqrt{\bar{r}S}$ is a finite $A/\sqrt{\bar{r}A}$-module via $\phi$. The ring $S/\sqrt{\bar{r}S}$ is catenary since $S$ is complete.\par 
Setting $B=A/\sqrt{\bar{r}A}$ and $C=S/\sqrt{\bar{r}S}$, we have that $B\subset C$  is an injective finite extension of equidimensional local rings and $C$ is catenary. The graph associated to $C$ is connected, and we need to show that the graph associated to $B$ is connected. This is shown below.\par
The graph associated to an equidimensional local ring is connected if and only if for every pair of minimal primes $P_{\alpha}$ and $P_{\beta}$ there is a sequence of prime ideals $P_1,\cdots,P_k$ such that, setting $P_{\alpha}=P_0$ and $P_{\beta}=P_{k+1}$, we have that $ht(P_i+P_{i+1})\leq 1$ (i.e. $ht(P_i+P_{i+1})=1$ if $P_i\neq P_{i+1}$) for all $0\leq i\leq k$. Accordingly, let $P_{\alpha}$ and $P_{\beta}$ be two arbitrary minimal prime ideals in $B$. Then we have prime ideals $\tilde{P}_{\alpha}$ and $\tilde{P}_{\beta}$ in $C$ lying over $P_{\alpha}$ and $P_{\beta}$, respectively, and $\tilde{P}_{\alpha}$ and $\tilde{P}_{\beta}$ are minimal in $C$ as well, by the going-up theorem. The graph associated to $C$ is connected, hence there exists a sequence of minimal prime ideals $\tilde{P_1},\cdots, \tilde{P_k}$ in $C$ so that, setting $\tilde{P}_0=\tilde{P}_{\alpha}$ and $\tilde{P}_{k+1}=\tilde{P}_{\beta}$,
$$ht_C(\tilde{P_i}+\tilde{P_{i+1}})=1,\ for\ 0\leq i\leq k.$$
Let $P_1,\cdots, P_k$ be the pullback of $\tilde{P_1},\cdots, \tilde{P_k}$ in $B$. We let $P_0=P_{\alpha}$ and $P_{k+1}=P_{\beta}$. To show that the graph associated to $B$ is connected, it is enough to show that $$ht_B(P_i+P_{i+1})\leq 1\ for\ all\ i\leq k.$$ This amounts to showing that if $P_i\neq P_{i+1}$ then $ht_B(P_i+P_{i+1})=1$. Accordingly, we assume that $P_i\neq P_{i+1}$. Since $ht_C(\tilde{P_i}+\tilde{P}_{i+1})=1$, there exists a prime ideal $\tilde{Q}$ in $C$ with height 1 containing $\tilde{P}_i$ and $\tilde{P}_{i+1}$. Let $Q$ be the pullback of $\tilde{Q}$ in $B$. Since $dim(B)=dim(C)=d-1$ and $C$ is catenary, equidimensional and $ht_C(\tilde{Q})=1$, we will have a chain of prime ideals
$$\tilde{Q}\subsetneq \tilde{Q_1}\subsetneq \cdots \subsetneq \tilde{Q}_{d-2}.$$Taking the pullback of this chain in $B$, we will have 
$$Q\subsetneq Q_1\subsetneq \cdots \subsetneq Q_{d-2}.$$
Hence, $ht_B(Q)\leq 1$, i.e. $ht_B(P_i+P_{i+1})\leq 1$. This shows that the graph associated to $B$ is connected and completes the proof of Proposition 2.2 and Theorem 1.4.
\end{proof}

In conclusion, we give an application of our results to projective schemes over a field. For any projective scheme $X$ of dimension $d$ over a field $k$, we can write $X$ as $Proj(k[x_0,\cdots,x_n]/I)$ for some $n$ and $I$ homogeneous in $k[x_0,\cdots,x_n]$, i.e. we have an embedding $X\hookrightarrow \mathbb{P}^n_k$. Let $A$ denote the local ring $(k[x_0,\cdots,x_n]/I)_{(x_0,\cdots,x_n)}$. Since $A$ is a local ring containing a field, we can consider the Lyubeznik numbers of $A$.\par

Our Main Theorem provides some supporting evidence for a positive answer to the open question whether the Lyubeznik numbers of the above ring $A$, $\lambda_{i,j}(A)$, depend only on the integers $i$, $j$ and the scheme $X$ but are independent of the embedding $X\hookrightarrow \mathbb{P}^n_k$ \cite[p.133]{l2}. Indeed, we have the following theorem which is a direct consequence of our Main Theorem.
\newtheorem*{2.7}{\bf Theorem 2.7}
\begin{2.7}
Let $X$ be an arbitrary projective scheme of dimension $d$. Under some embedding $\iota: X\hookrightarrow \mathbb{P}^n_k$, we can write $X=Proj(R)$, where $R=k[x_0,\cdots,x_n]/I$ with some homogeneous ideal $I$ in the polynomial ring $k[x_0,\cdots,x_n]$. Let $A:=R_{(x_0,\cdots,x_n)}$. Then $\lambda_{d+1,d+1}(A)$ does not depend on the choice of $n$ and $I$, i.e., it does not depend on the embedding $\iota: X\hookrightarrow \mathbb{P}^n_k$. In other words, it is a numerical invariant on $X$. Indeed, let $k^{sep}$ be the separable closure of $k$ and let $X_1,\cdots,X_s$ be the $d$-dimensional irreducible components of $X\times_kk^{sep}$. Let $\Gamma_X$ be the graph on vertices $X_0,\cdots,X_s$ and $X_i$, $X_j$ are joined by an edge if and only if $dim(X_i\cap X_j)<d-1$. Then $\lambda_{d+1,d+1}(A)$ equals the number of connected components of $\Gamma_X$. 
\end{2.7}

\end{document}